\documentclass{amsart}
\usepackage{graphics}

\newtheorem{theorem}{Theorem}[section]
\newtheorem{lemma}[theorem]{Lemma}

\theoremstyle{remark}

\numberwithin{equation}{section}

\DeclareMathSymbol{\nmid}{\mathbin}{AMSb}{"2D}

\begin{document}
\newcommand{\beqs}{\begin{equation*}}
\newcommand{\eeqs}{\end{equation*}}
\newcommand{\beq}{\begin{equation}}
\newcommand{\eeq}{\end{equation}}
\newcommand{\bin}[2]{\genfrac{(}{)}{0pt}{0}{#1}{#2}}
\newcommand\mymod[1]{\mbox{ mod}\ {#1}}
\newcommand\smymod[1]{\mbox{\scriptsize mod}\ {#1}}
\newcommand\mylabel[1]{\label{#1}}
\newcommand\eqn[1]{(\ref{eq:#1})}

\title[\tiny On representation of an integer as the sum of three squares ...]
{On representation of an integer as the sum of three squares and ternary quadratic forms with the discriminants $p^2$, $16p^2$}

\author{Alexander Berkovich}
\address{Department of Mathematics, University of Florida, Gainesville,
Florida 32611-8105}
\email{alexb@ufl.edu}
\thanks{The first author was supported in part by NSA grant H98230-09-1-0051}
\author{Will Jagy}
\address{Math.Sci.Res.Inst., 17 Gauss Way, Berkeley, CA 94720-5070}
\email{jagy@msri.org}

\date{January 30, 2011}


\begin{abstract}

Let $s(n)$ be the number of representations of $n$ as the sum of three squares.
We prove a remarkable new identity for $s(p^2n)- ps(n)$ with $p$ being an odd prime.
This identity makes nontrivial use of ternary quadratic forms with discriminants $p^2$, $16p^2$.
These forms are related by Watson's transformations. To prove this identity we employ the Siegel--Weil 
and the Smith--Minkowski product formulas.

\end{abstract}
\dedicatory{} 

\maketitle

\section{Introduction} \label{sec:intro}
\medskip

Let $(a,b,c,d,e,f)(n)$ denote the number of integral representations of $n$ by the positive ternary quadratic form
$ax^2+by^2+cz^2+dyz+ezx+fxy$. We will take $(a,b,c,d,e,f)(n)=0$, whenever $n \not\in\mathbf  N$.
The discriminant $\Delta$ of a ternary form $ax^2+by^2+cz^2+dyz+ezx+fxy$ is defined as 
\beqs
\Delta=\frac{1}{2}\det
\begin{bmatrix} 2a & f & e \\ f & 2b & d \\ e & d & 2c
\end{bmatrix}
= 4abc +def -ad^2-be^2-cf^2.
\eeqs 
We say that two ternary quadratic forms $\tilde{f}(x,y,z)$ and $\tilde{g}(x,y,z)$ with the discriminant  $\Delta$  are in the same genus 
if they are equivalent over $\mathbf Q$ via a matrix in $GL(3,\mathbf Q)$ whose entries have denominators prime to $2\Delta$.
We add that  this is the case if and only if these  forms  are equivalent over the real numbers and over the $p$-adic integers  $\mathbf Z_{p}$
for all primes $p$ \cite{Cas}, \cite{Jones}, \cite{Wat}.\\

It is well known that all ternary forms with discriminant $4$ are equivalent to $x^2+y^2+z^2$ \cite{Dick}, \cite{Landau}.
Let $p$ be an odd prime. Lehman derived elegant counting formulas for ternary genera  in \cite{Leh}. Using his results, 
it is straightforward to check that all ternary forms with the discriminant $p^2$ belong to the same genus, say $TG_{1,p}$.
There are twelve genera of  ternary forms with the discriminant $16p^2$. However, if one imposes additional constraints
on the forms with  $\Delta=16p^2$, namely
\beqs
(a,b,c,d,e,f)(n) = 0, \mbox{   when } n\equiv 1,2\pmod 4,
\eeqs
\beqs
d\equiv e\equiv f\equiv 0 \pmod 2,
\eeqs
then we will show in Section 6 that all these ternaries belong to the same genus, say $TG_{2,p}$.
In Section 8 we will show how to relate  $TG_{1,p}$ and  $TG_{2,p}$ and Watson's transformation.

Let $s(n)$ denote the number of representations of $n$ by ternary form  $x^2+y^2+z^2$, so
\beqs
s(n) = (1,1,1,0,0,0)(n).
\eeqs
In \cite{Brk} the first author utilized $q$-series techniques to prove the following two theorems: 
\begin{theorem}\label{nt1}
\beq
s(9n)-3s(n) = 2(1,1,3,0,0,1)(n)-4(4,3,4,0,4,0)(n).
\mylabel{eq:1.1}
\eeq
\end{theorem}
\noindent
\begin{theorem}\label{nt2}
\beq
s(25n)-5s(n) =  4(2,2,2,-1,1,1)(n)-8(7,8,8,-4,8,8)(n).
\mylabel{eq:1.2}
\eeq
\end{theorem}
\noindent
Our main object here is to prove the following 
\begin{theorem}\label{nt3}
Let $p$ be an odd prime, then
\beq
s(p^2n)-p s(n) =  48\sum_{\tilde{f}\in TG_{1,p}}\frac{R_{\tilde{f}}(n)}{|\mbox{Aut}(\tilde{f})|} -96\sum_{\tilde{f}\in TG_{2,p}}\frac{R_{\tilde{f}}(n)}{|\mbox{Aut}(\tilde{f})|},
\mylabel{eq:1.3}
\eeq
where ${|\mbox{Aut}(\tilde{f})|}$ denotes the number of integral automorphs of a ternary form $\tilde{f}$, $R_{\tilde{f}}(n)$ denotes the number of representations of $n$ by $\tilde{f}$,
and a sum over forms in a genus should be understood to be the finite sum resulting from taking a single representative from each equivalence class of forms.
\end{theorem} 
\noindent
This theorem was first stated in \cite{Brk}. We remark that, somewhat similar in flavor, the so-called $S$-genus identities were recently discussed in \cite{BeJ}, \cite{BHJ}.
In what follows we will require the following 
\begin{theorem}\label{nt4}
\beq
s(n) =\frac{16}{\pi}\sqrt n\psi(n) L(1,\chi(n))P(n),
\mylabel{eq:1.4}
\eeq
where for $n = 4^a k$, $4 \nmid k$ one has 
\beq
\psi(n)=
\begin{cases}
0 & \text{if }k\equiv 7 \pmod 8,\\
2^{-a} & \text{if } k\equiv 3\pmod 8,\\
3\cdot 2^{-a-1} & \text{if } k\equiv 1,2 \pmod 4;
\end{cases}
\mylabel{eq:1.4a}
\eeq
$L(1,\chi(n))=\sum_{m=1}^\infty(-4n|m)m^{-1}$ with $\chi(n)=(-4n|\bullet)$, 
the Kronecker symbol and
\beq
P(n)=\prod_{(p'^{2})^{b}||n}\biggl(1+\frac{1}{p'}+\frac{1}{p'^2}+\cdots +\frac{1}{p'^{b-1}}+\frac{1}{p'^b\bigl(1-\bigl(
-np'^{-2b}|p'\bigr)p'^{-1}\bigr)}\biggr),
\mylabel{eq:1.4b}
\eeq 
with the product over all odd primes $p'$ such that  $p'^{2}\mid n$.
\end{theorem}
\noindent
Proofs of this Theorem may be found in \cite{Bach} and \cite{PB}.
We observe that $L(1,\chi(n))$ can be written as the infinite product 
\beq
L(1,\chi(n))= \frac{\pi^2}{8}\prod_{p'}\biggl(1+\bigl(-n|p'\bigr)\frac{1}{p'}\biggr),
\mylabel{eq:1.5}
\eeq
where $p'$ runs through all odd primes.\\
Before we move on we comment  that for squarefree $n\equiv 3 \pmod 8$, $n\ge 11$
\beqs
L(1,\chi(n))= \frac{3}{2} \frac{\pi}{\sqrt{n}}h(n),
\eeqs
where $h(n)$ is the class number of the quadratic field $Q(\sqrt{-n})$.

\bigskip
\section{The Siegel--Weil formula for ternary quadratic forms} \label{sec:2}
\medskip
 
Let $T$ be a genus of positive ternary forms with the discriminant $\Delta$.
Then the Siegel--Weil formula \cite{Siegel} implies that 
\beq
\sum_{\tilde{f}\in T}\frac{R_{\tilde{f}}(n)}{|\mbox{Aut}(\tilde{f})|}= 4 \pi M(T)\sqrt\frac{n}{\Delta}\prod_{p'}d_{T,p'}(n),
\mylabel{eq:2.1}
\eeq
where ${|\mbox{Aut}(\tilde{f})|}$ denotes the number of integral automorphs of a ternary form $\tilde{f}=ax^2+by^2+cz^2+dyz+ezx+fxy$, 
while $R_{\tilde{f}}(n)$ denotes the number of representations of $n$ by $\tilde{f}$.
The sum on the left is over forms in a genus. Again, this sum (here and everywhere) should be interpreted as the finite sum resulting 
from taking a single representative from each equivalence class of forms. 
The product on the right is over all primes, the mass of the genus is defined by
\beqs
M(T):=\sum_{\tilde{f}\in T}\frac{1}{|\mbox{Aut}(\tilde{f})|},
\eeqs
and $d_{T,p'}(n)$ denotes the $p'$-adic (local) representation density, defined by 
\beqs
d_{T,p'}(n):=  \frac{1}{p'^{2t}}|\mbox\{(x,y,z)\in Z^{3} : ax^2+by^2+cz^2+dyz+ezx+fxy \equiv n \pmod {p'^t}\}|,
\eeqs
for  sufficiently large $t$. We comment that  $ax^2+by^2+cz^2+dyz+ezx+fxy$  can be chosen to be any form $\in T$.
In \cite{Siegel} Siegel proved that when $\gcd(2\Delta,p')=1$
\beq
d_{T,p'}(n) =    
\begin{cases}
\big(\frac{1}{p'}+1\big)+\frac{1}{p'^{k+1}}((-m|p')-1) & \text{if } n= m p'^{2k}, \;\;  p' \nmid m,\\
\big(\frac{1}{p'}+1\big)\big(1-\frac{1}{p'^{k+1}}\big) & \text{if } n= mp'^{2k+1}, \;\; p' \nmid m.
\end{cases}
\mylabel{eq:2.2}
\eeq
It is not hard to check that \eqn{1.4} follows easily from \eqn{2.1} and  \eqn{2.2}, provided one recognizes that 
\beq
\psi(n)= d_{x^2+y^2+z^2,2}(n),
\mylabel{eq:2.3}
\eeq
where $\psi(n)$ is defined in \eqn{1.4a}.
It is easy to check that 
\beqs
d_{x^2+y^2+z^2,2}(4n) = \frac{1}{2}d_{x^2+y^2+z^2,2}(n),
\eeqs
and that 
\beq
d_{x^2+y^2+z^2,2}(n)= 0 ,\;\;\; \text{if }n\equiv 7\pmod 8.
\mylabel{eq:2.4}
\eeq
It remains to verify that 
\beq
d_{x^2+y^2+z^2,2}(n)=
\begin{cases}
1 & \text{if } n\equiv 3\pmod 8,\\
\frac{3}{2} & \text{if } n\equiv 1,2 \pmod 4.
\end{cases}
\mylabel{eq:2.5}
\eeq
This can be easily accomplished with the  help of the following 
\begin{lemma}\label{l1}
The number of roots of 
\beqs
x^2 \equiv c\pmod {2^t},\;\;\; 3\le t,\;\;\; c\equiv 1 \pmod 2
\eeqs
is four or zero, according as $c\equiv 1\pmod 8$ or $c\not\equiv 1\pmod 8$. 
\end{lemma}
\noindent
Proof of this Lemma may be found in \cite{Dick}.
Next, we observe that 
\begin{align*}
& \text{ }|\{(x,y,z)\in Z^{3} : 0\leq x,y,z<2^t, x^2+y^2+z^2 \equiv n \pmod {2^t}\}| = \\ 
& 4 |\{(y,z)\in Z^{2} : 0<y,z<2^{t}, yz\equiv 1\pmod 2 \}|=4\cdot 2^{t-1}\cdot 2^{t-1},
\end{align*}
when $n\equiv 3\pmod 8$.
Hence,
\beqs
d_{x^2+y^2+z^2,2}(n)= \frac{4}{2^{2t}} 2^{t-1} 2^{t-1} =1,
\eeqs
when $n\equiv 3\pmod 8$.
Analogously, when $n\equiv 1\pmod 8$, $3\le t$ we find that 
\begin{align*}
& \text{ }|\{(x,y,z)\in Z^{3} : 0\leq x,y,z<2^{t}, x^2+y^2+z^2 \equiv n\pmod{2^t} \}|   \\
=&        3|\{(x,y,z)\in Z^{3} : 0\leq x,y,z<2^{t}, x\equiv 1\pmod 2, x^2+y^2+z^2 \equiv n\pmod{2^t}\}|  \\ 
=&        3\cdot4 |\{(y,z)\in Z^{2}:0\leq y,z<2^{t}, yz\equiv 0\pmod 2, y\equiv z\pmod 4 \}| \\
=&        3\cdot4\cdot\frac{1}{2}\cdot2^{t-1}\cdot 2^{t-1}. 
\end{align*}
And so 
\beqs
d_{x^2+y^2+z^2,2}(n)= \frac{1}{2^{2t}}\cdot3\cdot4\cdot\frac{1}{2}\cdot2^{t-1}\cdot 2^{t-1} = \frac{3}{2} \;\;\; \text{when } n\equiv 1\pmod 8,
\eeqs
as desired. Other cases in \eqn{2.5} can be handled in a very similar manner.\\
We can now rewrite \eqn{1.4} as 
\beqs
s(n)= 2\pi \sqrt{n} \prod_{p'}d_{x^2+y^2+z^2,p'}(n).
\eeqs
Consequently, 
\beq
s(p^2n)-p s(n) = 2\pi \sqrt{n} \psi(n)\Gamma_{p}(n)\prod_{\gcd(p',2p)=1}d_{x^2+y^2+z^2,p'}(n),
\mylabel{eq:2.6}
\eeq
where 
\beqs
\Gamma_{p}(n): = p(d_{x^2+y^2+z^2,p}(p^2n)-d_{x^2+y^2+z^2,p}(n)).
\eeqs
From \eqn{2.2} we have at once 
\beq
\Gamma_{p}(n)=
\begin{cases}
\frac{p-1}{p^{1+k}}(1-(-m|p)) & \text{if } n= m p^{2k}, \; \; p \nmid m, \\
\frac{p-1}{p^{1+k}}\big(1+\frac{1}{p}\big) & \text{if } n= mp^{2k+1}, \;\; p \nmid m.
\end{cases}
\mylabel{eq:2.7}
\eeq

\bigskip
\section{ Computing some local representation densities.\\ The non-dyadic case.} \label{sec:3}
\medskip

In this section we prove 
\begin{theorem}\label{nt5}
Let $p$ be an odd prime and $u$ be any integer with $(-u|p)=-1$. Let $G$ be some ternary 
genus such that  $$ \tilde{f} \sim_p ux^2+ p(y^2+ u z^2) $$ for any $\tilde{f}$ in $G$. Then
\beq
d_{G,p}(n)=  \frac{p}{p-1} \Gamma_{p}(n).
\mylabel{eq:3.1}
\eeq
Here, and everywhere, the relation  $\tilde{f} \sim_p \tilde{g}$ means that the two quadratic forms $\tilde{f}$ and $\tilde{g}$ are 
equivalent over the $p$-adic integers $\mathbf Z_{p}$.
\end{theorem}
\noindent
It may not be obvious that the above theorem is identical to the special case $ \epsilon_{p} = -1$  of  Lemma 4.2 in~\cite{BHJ}. 
Here we take a leisurely and self-contained approach, counting solutions of the relevant equation modulo $p^t$ for large $t$.
Suppose 
\beqs
ux^2+ p(y^2+ u z^2)\equiv\ n \pmod{p^t}\text{ with } p^2 \mid n,\;\;\;2\le t.
\eeqs
Then, thanks to $(-u|p)=-1$, we have $p \mid x$, $p \mid y$, $p \mid z$.  
This observation implies that 
\beq
d_{G,p}(np^{2k})= \frac{d_{G,p}(n)}{p^{k}},\;\; \text{if } p^2\nmid n.
\mylabel{eq:3.2}
\eeq
Hence \eqn{3.1} holds true for all $n$ if it holds true for all $n$ such that $p^2\nmid n$.
There are two cases to consider.
First, when $p\nmid n$ we have that 
\begin{align*}
&\text{ }|\{ (x,y,z)\in Z^{3} : 0\le\ x,y,z< p^{t},  ux^2+ p(y^2+ u z^2) \equiv n \pmod {p^t}\}|  \\ 
=&|\{(x,y,z)\in Z^{3} :0\le\ x,y,z< p^{t}, x^2 \equiv un -up(y^2+ u z^2) \pmod {p^t}\}| \\ 
=& \sum_{y=0}^{p^{t}-1}\sum_{z=0}^{p^{t}-1} ( 1+ ((un -up(y^2+ u z^2))|p))= p^{2t}( 1+ (un|p))=p^{2t}( 1-(-n|p)).
\end{align*}
And so 
\beq
d_{G,p}= \frac{1}{p^{2t}} p^{2t}( 1-(-n|p))= ( 1-(-n|p)).
\mylabel{eq:3.3}
\eeq
We comment that in the discussion  above  we used a well known 
\begin{lemma}\label{l2}
Let $p'$ be an odd prime not dividing $c$.
The number of roots of 
\beqs
x^2 \equiv c\pmod {p'^t},\;\;\; t\ge 1 
\eeqs
is the same as the number $(0$ or $2)$ of roots when $t=1$. That is 
\beqs
|\{0\le\ x<p^{t}: x^2 \equiv c \pmod {p^t}\}|=1 + (c|p).
\eeqs
\end{lemma}
\noindent
This lemma is proven in \cite{Dick}.\\
Second, when $n=pm$ and $p\nmid m$ we have that 
\begin{align*}
& \text{ }|\{ (x,y,z)\in Z^{3} : 0\le\ x,y,z< p^{t},  ux^2+ p(y^2+ u z^2) \equiv pm \pmod {p^t}\}|\\
=& p^2\text{ }|\{(x,y,z)\in Z^{3} : 0\le\ x,y,z< p^{t-1}, upx^2+ y^2+ u z^2 \equiv m \pmod {p^{t-1}}\}|\\
=& p^2 |\{(x,y,z)\in Z^{3} :0\le\ x,y,z< p^{t-1}, y^2 \equiv m -upx^2-u z^2 \pmod {p^{t-1}}\}|\\
=& p^2 \sum_{x=0}^{p^{t-1}-1}\sum_{z=0}^{p^{t-1}-1} ( 1+ ( (m -upx^2-u z^2) |p))\\
=& p^{2t} + p^2\sum_{x=0}^{p^{t-1}-1}\sum_{z=0}^{p^{t-1}-1} ( (m -upx^2-u z^2)|p)\\
=& p^{2t} + p^{t+1}\sum_{z=0}^{p^{t-1}-1}( (m -u z^2)  |p)= p^{2t}- p^{t+1}\sum_{z=0}^{p^{t-1}-1}((-um +z^2) |p)\\
=& p^{2t}+ p^{t+1} p^{t-2} = p^{2t}\big(1+\frac{1}{p}\big).
\end{align*}
This time we used another well known fact:
\beq
\sum_{y=0}^{p-1}((y^2 + a)|p)= -1,
\mylabel{eq:3.4}
\eeq
with $p$ being an odd prime not dividing $a$. \\
Hence
\beq
d_{G,p}= \frac{1}{p^{2t}} p^{2t}\big(1+\frac{1}{p}\big) =\big(1+\frac{1}{p}\big),
\mylabel{eq:3.5}
\eeq
as desired. \\
Our proof of Theorem 3.1 is now complete.

\bigskip
\section{ Computing some local representation densities.\\ The dyadic case.} \label{sec:4}
\medskip

In this section we prove two theorems.
\begin{theorem}\label{nt6}
Let $G_{1}$ be some ternary genus such that  
$$
\tilde{f} \sim_2 yz-x^2 
$$ 
for any $\tilde{f}$ in $G_{1}$.
Let $n = 4^a k, 4 \nmid k$, then
\beq
d_{G_{1},2}(n)=
\begin{cases}
\frac{3}{2} & \text{if }k\equiv 7\pmod 8,\\
\frac{3}{2}-\frac{1}{2^{a+1}} & \text{if } k\equiv 3\pmod 8,\\
\frac{3}{2}-\frac{3}{2^{a+2}}& \text{if } k\equiv 1,2 \pmod 4.
\end{cases}
\mylabel{eq:4.1}
\eeq
\end{theorem}
\noindent
\begin{theorem}\label{nt7}
Let $G_{2}$ be some ternary 
genus such that  $$ \tilde{f} \sim_2 4yz -x^2 $$ for any $\tilde{f}$ in $G_{2}$.
Let  $n = 4^a k, 4 \nmid k$, then
\beq
d_{G_{2},2}(n)=
\begin{cases}
3& \text{if }k\equiv 7\pmod 8,\\
3-\frac{1}{2^{a-1}} & \text{if } k\equiv 3\pmod 8,\\
3-\frac{3}{2^{a}}& \text{if } k\equiv 1,2 \pmod 4.
\end{cases}
\mylabel{eq:4.2}
\eeq
\end{theorem}
\noindent
Comparing \eqn{1.4a}, \eqn{4.1}, and \eqn{4.2} we have at once
\beq
\psi(n)= 2 d_{G_{1},2}(n) - d_{G_{2},2}(n).
\mylabel{eq:4.3}
\eeq
We note two related recurrences
\beq
2 d_{G_1,2}(n) -  d_{G_2,2}(4 n) = 0
\mylabel{eq:4.4}
\eeq
and
\beq
4 d_{G_1,2}(n) -  d_{G_2,2}(n) = 3.
\mylabel{eq:4.5}
\eeq
To prove Theorem 4.1 and Theorem 4.2, it is sufficient to  show  that \eqn{4.4} and \eqn{4.5},
together with the initial conditions 
\beq
d_{G_{2},2}(n)=
\begin{cases}
3& \text{if }n\equiv 7\pmod 8,\\
1 & \text{if } n\equiv 3\pmod 8,\\
0& \text{if } n\equiv 1,2 \pmod 4,
\end{cases}
\mylabel{eq:4.6}
\eeq
hold  true.
Note that \eqn{4.4} follows easily from 
\begin{align*}
& \text{ }|\{(x,y,z)\in Z^{3}: 0\le x,y,z<2^{t},4yz-x^2 \equiv 4n \pmod {2^t}\}| \\ 
=& 2\cdot 4\cdot 4 |\{(x,y,z)\in Z^{3}:0\le x,y,z<2^{t-2},yz-x^2\equiv n \pmod {2^{t-2}}\}|.
\end{align*}
Clearly, when $n\equiv 1,2 \pmod 4$ we have 
\beqs
4yz-n\equiv 2,3,6,7 \pmod 8.
\eeqs
Recalling Lemma 2.1, we see that the congruence 
\beqs
4yz-x^2 \equiv n \pmod {2^t} 
\eeqs
has no solutions when $t\ge3$.
Consequently,
\beqs
d_{G_{2},2}(n)= 0 \;\;\; \text{if } n\equiv 1,2 \pmod 4.
\eeqs
Next, when $n\equiv 3 \pmod 8$ we have 
\begin{align*}
& \text{  }|\{(x,y,z)\in Z^{3}:0\le x,y,z<2^{t},4yz-x^2 \equiv n \pmod {2^t}\}| \\ 
=& 4|\{(y,z)\in Z^{2}: 0\le y,z<2^{t},yz\equiv 1 \pmod 2 \}|= 4\cdot 2^{t-1} 2^{t-1} = 2^{2t}.
\end{align*}
Hence
\beqs
d_{G_{2},2}(n)= 1 \;\; \text{if } n\equiv 3 \pmod 8.
\eeqs
The case $n\equiv 7 \pmod 8$ in \eqn{4.6}  can be treated in an analogous manner.
\begin{align*}
& \text{ }|\{(x,y,z)\in Z^{3}:0\le x,y,z<2^{t},4yz-x^2 \equiv n \pmod {2^t}\}|  \\ 
=& 4\; |\{(y,z)\in Z^{2} : 0\le y,z<2^{t}, yz\equiv 0 \pmod 2\}|= 4\cdot(2^{t} 2^{t}- 2^{t-1} 2^{t-1} )= 3\cdot 2^{2t}.
\end{align*}
Hence
\beqs
d_{G_{2},2}(n)= 3 \; \text{if } n\equiv 7 \pmod 8.
\eeqs
And so we have established the initial conditions \eqn{4.6}. It remains to prove \eqn{4.5}.
We shall require the following easy companion to Lemma 2.1
\begin{lemma}
Let $t = 2 s + 3 + \delta$, with integers $0 \leq \delta \leq 1$ and $0\le s$. \\ 
Let $S_{t}(c): = |\left\{0 \leq x < 2^t: \;\; x^2 \equiv c \pmod {2^t} \right\}|$. \\
Then
$$
S_{t}(4^m \cdot ( 8 r + 1))=\left\{
\begin{array}{ll}
2^{m+2} & \mathrm{if } \; \; 0 \leq m \leq s ,\\
2^{s+1 + \delta} & \mathrm{if } \;\; m = s + 1, \\ 
2^{s+1 + \delta} & \mathrm{if } \;\; m = s + 2.\\
\end{array}
\right.
$$
When $m = s+2$ the formula $4^m \cdot ( 8 r + 1)$
refers to $0$, as then $2 s + 4 \geq t$ and $4^m \geq 2^t$; in fact $2^t\mid 4^m$. 
It is important to note that no $0 \leq c < 2^t$ other than the specified values above are allowed to have $S_{t}(c)\neq 0$.
\end{lemma}
\noindent
To proceed further we define 
\beqs
P_{i,t}(c) = |\left\{ (y,z)\in Z^{2} :0 \leq y,z < 2^t,4^{i-1} y z \equiv c \pmod {2^t} \right\}|,\;\;\;i=1,2,
\eeqs
and 
\beqs
C_{m,t} = \left\{ c\in Z : 0 \leq c < 2^t,c \equiv 4^m \pmod {8 \cdot 4^m} \right\},\;\;\; m\in Z.
\eeqs
Again we comment that when $ m = s+2$ the condition $c \equiv 4^m \pmod {8 \cdot 4^m}$ means $c = 0$.
It is not hard to check that 
\beq
4 \cdot P_{1,t}( n) - P_{2,t}(n)=
\begin{cases}
2^{t+1} & \text{if } n\equiv 1\pmod 2,\\
2^{t+2} & \text{if } n\equiv 0\pmod 2,
\end{cases}
\mylabel{eq:4.7}
\eeq
and that for $t = 2 s + 3 + \delta$, with integers $0 \leq \delta \leq 1$ and $0\le s$
\beq
|C_{m,t}|=
\begin{cases}
2^{2s-2m + \delta}& \text{if } 0 \leq m \leq s, \\
1 & \text{if }  m = s + 1 ,\\
1 & \text{if }  m = s + 2 .
\end{cases}
\mylabel{eq:4.8}
\eeq
Next, we define for $i=1,2$
\beqs
L_{i,t}(c) = |\{(x,y,z)\in Z^{3}: 0 \leq x,y,z < 2^t,\;\;\; 4^{i-1}y z - x^2 \equiv c \pmod {2^t}\}|.
\eeqs
From Lemma 4.3 it is easy to see that 
\beqs
L_{i,t}(n) = \sum_{m=0}^{s+2}\sum_{c \in C_{m,t}}S(c) P_{i,t}(n+c),\;\;\; i=1,2.
\eeqs
Making use of Lemma 4.3, \eqn{4.7}, and \eqn{4.8} we find that 
\beq
4 \cdot L_{1,t}(n) - L_{2,t}(n) = \sum_{m=0}^{s+2}\sum_{c \in C_{m,t}}S(c)\cdot(4\cdot P_{1,t}(n+c) - P_{2,t}(n+c))= 3\cdot 4^t.
\mylabel{eq:4.9}
\eeq
Finally, we note that for sufficiently large $t$
\beq
d_{G_{i},2}(n)=  \frac{1}{4^t} L_{i,t}(n),\;\;\; i=1,2.
\mylabel{eq:4.10}
\eeq
Combining \eqn{4.9} and \eqn{4.10} we see that \eqn{4.5} holds true.
Our proofs of Theorem 4.1 and Theorem 4.2 are now complete.

\bigskip
\section{Computing the mass of the ternary genus $TG_{1,p}$} \label{sec:5}
\medskip

In this section we prove that
\beq
M(TG_{1,p})= \frac{p-1}{48},
\mylabel{eq:5.1}
\eeq
where $p$ is a fixed odd prime and 
\beqs
M(TG_{1,p}):=\sum_{\tilde{f}\in TG_{1,p}}\frac{1}{|\mbox{Aut}(\tilde{f})|}.
\eeqs
To prove \eqn{5.1} we will employ the Smith–-Minkowski-–Siegel mass formula.
This formula gives the mass as an infinite product over all primes.
Many published versions of this  formula have small errors.
In this paper we will follow a reliable account by Conway and Sloane \cite{CS}. 
From equation (2) in \cite{CS} we have that  
\beq
M(TG_{1,p})=  \frac{1}{\pi^2} \prod_{p'} 2m_{p'},
\mylabel{eq:5.2}
\eeq
where $p'$ runs through all primes and where local masses $m_p'$ are
defined in equation (3) in \cite{CS} by 
\beq
m_{p'}= \prod_{q} M_q  \prod_{q<q'} {(q'/q)}^{\frac{n(q)n(q')}{2}} 2^ { n(I,I)-n(II)}.
\mylabel{eq:5.3}
\eeq
Here $q$ ranges over all  powers $p'^t$ of $p'$ (including those with negative $t$).
The last factor in \eqn{5.3} is  $1$  for all odd primes.
So if the $p'$-adic Jordan decomposition of $\tilde{f}\in TG_{1,p}$ is given by 
\beqs
\sum_{q} q \tilde{f}_q,
\eeqs
then 
\beqs
n(q)= \dim(\tilde{f}_q).
\eeqs
For all odd primes $p'$ such that $p'\nmid p$, the $p'$-adic Jordan decomposition of 
any form $\tilde{f}\in TG_{1,p}$ can be taken to be  $(x^2+y^2+z^2)$; this follows from Theorem 29 in \cite{Wat}.
So with the aid of Table 2 in \cite{CS} we find that 
\beqs
n(1)=3,\;\;\;\; M_1= \frac{p'^2}{2(p'^2-1)}.
\eeqs
If $q\neq 1$ then 
\beqs
n(q)=0,\;\;\;\; M_q= 1.
\eeqs
Hence 
\beq
m_{p'} = M_1 =  \frac{p'^2}{2(p'^2-1)},\;\; p' \nmid 2p.
\mylabel{eq:5.4}
\eeq
Next, the $p$-adic Jordan decomposition of any form  $\tilde{f}\in TG_{1,p}$ is given by
$\tilde{f}_1+p\tilde{f}_p$ with $\tilde{f}_1=u x^2$, $\tilde{f}_p=(y^2+ uz^2)$,
for a unit $u$ satisfying $(-u|p)=-1$ (see Section 6).\\
And so from Table 1 and Table 2 in \cite{CS}, we see that 
$n(1)=1$, species$(1)=1$, $M_1=\frac{1}{2}$, and 
$n(p)=2$, species$(p)=2$, $M_p=\frac{p}{2(p+1)}$. 
And so we find that 
\beq
m_p = M_1 M_p (p/1)^\frac{2}{2}  =\frac{p^2}{4(p+1)}.
\mylabel{eq:5.5}
\eeq
Finally, one possible $2$-adic Jordan decomposition of any form $\tilde{f}\in TG_{1,p}$ is  given  by 
\beqs
\frac{1}{2}\tilde{f}_{\frac{1}{2}}+\tilde{f}_1+2\tilde{f}_2, 
\eeqs
with $\tilde{f}_{\frac{1}{2}}=2yz$, $\tilde{f}_1 = -x^2$, $\tilde{f}_2 = 0$. 
This follows from Theorem 29 in \cite{Wat}.
We note that $\tilde{f}_2 $ is a bound love form. It contributes a factor of $\frac{1}{2}$ to the mass
\beqs
M_{2}= \frac{1}{2}.
\eeqs
Obviously $n\big(\frac{1}{2}\big)=2$, $n(1)=1$, and $n(2)=0$.\\
Next, 
$\tilde{f}_{\frac{1}{2}}$ is of the type II$_{2}$. It is bound and has octane value $=0$, species $=3$, $M_{\frac{1}{2}}=\frac{2}{3}$.
Also, $\tilde{f}_{1}$ is of the type I$_{1}$. It is free and has octane value $=0-1=-1$, species $=0+$, $M_{1}=1$.
In \eqn{5.3}, $n($I,I$)$ is the total number of pairs of adjacent constituents $\tilde{f}_{q}$, $\tilde{f}_{2q}$ that are both of type I
and $n($II$)$ is the sum of the dimensions of all Jordan constituents that have type II.
Clearly $n($I,I$)=0$ and $n($II$)=2$.
So 
\beq
m_{2} =  \frac{2}{3}(1) \frac{1}{2} (2/1)^{\frac{2}{2}} 2 ^{0-2}=\frac{1}{6}.
\mylabel{eq:5.6}
\eeq
Combining \eqn{5.2}--\eqn{5.6}, we obtain 
\beqs
M(TG_{1,p})=\frac{1}{\pi^2}\frac{1}{3}\frac{p^2}{2(p+1)}\prod_{\gcd(2p,p')=1}\frac{p'^2}{p'^2-1} =\frac{p-1}{8\pi^2}\prod_{p'}\frac{p'^2}{p'^2-1}. 
\eeqs
Recalling that
\beqs
\prod_{p'} \frac{p'^2}{p'^2-1}=\sum_{n\geq1}\frac{1}{n^2}=\frac{\pi^2}{6},
\eeqs
we see that \eqn{5.1} holds true.

\bigskip
\section{A Tale of Two Genera} \label{sec:6}
\medskip

Here we will give an overview of the construction of $TG_{2,p}$. We are given $TG_{1,p}$.
The sextuple 
$$ 
\langle a,b,c,d,e,f \rangle
$$ 
refers to 
$$
\hspace{4mm} a x^2 + b y^2 + c z^2 + d y z + e z x + f x y ,
$$ 
with Gram matrix
\begin{equation}
\left(  \begin{array}{ccc}
  2a & f & e\\
  f & 2b & d\\
  e & d & 2c  
\end{array} 
  \right). 
\end{equation}
First we will show that any form in $TG_{1,p}$ is equivalent to a form in Convenient Shape 1, which is just
$ \langle a,b,c,d,e,f \rangle$ with $a \equiv 3 \pmod 4$, then $d$ odd and $e,f$ even. Any primitive form represents an odd number, 
therefore it primitively represents an odd number $a$, so we may insist that $a$ be odd.
A particularly simple operation taking a form to an equivalent one is constructing the form with Gram matrix $M_{ij}' G M_{ij}$, 
where $G$ is the current Gram matrix of the form and $M_{ij}$ is the result of beginning with the identity matrix and placing a single $1$ at position $ij$, 
and $M_{ij}'$ denotes the transpose of  $M_{ij}$. We can also permute variables with the matrix
\begin{equation}
M_0 = \left(  \begin{array}{ccc}
  1 & 0 & 0\\
  0 & 0 & -1\\
  0 & 1 & 0  
\end{array} 
  \right). 
\end{equation}
If $e$ and $f$ are both even we are done. Otherwise, at least one of them is odd. If $f$ is the odd one apply $M_0,$ as in  $M_0' G M_0$, 
to arrive at $\langle a,c,b,-d,-f,e \rangle$; call these the new values of all the letters $a,b,c,d,e,f.$ 
If $d$ is even, apply $M_{12}$ to get odd $d$ in  $\langle a, a+b+f, c, d+e, e, f + 2a \rangle$. 
If $f$ is odd apply $M_{32}$ to get even $f$ in  $\langle a, b+c+d, c, d+2c, e, f + e \rangle$.
From the definition of the discriminant, $\Delta = 4 a b c + d e f - a d^2 - b e^2 - c f^2$, with
$\Delta =  p^2 \equiv 1 \pmod 4$; it follows that $b$ is even, so we have $b,f$ even and $a,d$ odd. 
Finally, apply $M_{21}$ to get even $e$ in  $\langle a+b+f, b, c, d, e + d, f + 2b \rangle$, where the new value of $a$ is still odd, 
while $e$ has become even. With $a,d$ odd and $e,f$ even, all the terms in  $\Delta = 4 a b c + d e f - a d^2 - b e^2 - c f^2$ are divisible by 4 except $-a d^2$. 
Since $d^2 \equiv 1 \pmod 4$ and $ \Delta \equiv 1 \pmod 4$, it follows that $ a \equiv 3 \pmod 4$.
Given a primitive form  $ \langle a,b,c,d, e,f \rangle$ in Convenient Shape 1,
define a mapping $ \Phi$ giving another primitive form by
$$ 
\Phi \,(\langle a,b,c,d, e,f \rangle) \, = \, \langle a,4b,4c,4d,2 e,2f \rangle.
$$
Note that if $ g(x,y,z) $ is in Convenient Shape 1 and $h = \Phi(g)$, then $ h(x,y,z) = g(x,2y,2z)$.
Any primitive form  $ \langle a,b,c,d,e,f \rangle$ with  $d,e,f \equiv 0 \pmod {2}$ that does not represent any number  
$n \equiv 1,2 \pmod {4}$  has $\Delta \equiv 0 \pmod {16}$ can be put in (is equivalent to a form in) Convenient Shape 2, 
that is with $b,c,d,e,f$ all divisible by $4$, and with  $ a \equiv 3 \pmod 4$. \\
To save space, we will introduce matrices
$$ E=
\left(  \begin{array}{rrr}
  1 & 0 & 0\\
  0 & 2 & 0\\
  0 & 0 & 2  
\end{array} 
  \right),  
$$
and
$$ D=
\left(  \begin{array}{rrr}
  2 & 0 & 0\\
  0 & 1 & 0\\
  0 & 0 & 1  
\end{array} 
  \right).  
$$
Note that $DE = ED = 2 I$.
At this point we will give an outline of the proof that  $\Phi$ is a well-defined bijection between  $TG_{1,p}$ and $TG_{2,p}$, that $\Phi$ preserves automorphs, 
that $TG_{2,p}$ is in fact a genus, and finally give the $p$-adic diagonalization for all these forms, along with a $2$-adic Jordan decomposition.
Given a primitive form $g$ in Convenient Shape 1, with Gram matrix $G$, the Gram matrix for $\Phi(g)$ is $EGE$. 
Furthermore $\Phi(g)$ is also primitive.
Now, suppose that forms $g,h$ in $TG_{1,p}$ are equivalent and are already in Convenient Shape 1. Suppose they have Gram matrices $G$ and $H$,respectively. So we are saying there is an integral matrix $P$ such that
$P' G P = H$.  It turns out that matrix elements $p_{21}, p_{31}$ are even, and so
$ \frac{1}{2} D P E$ is integral. But then
\begin{align}
\big(\frac{1}{2} D P E\big)' \big(EGE\big)  \big(\frac{1}{2} D P E  \big) & = \big(\frac{1}{2} E P' D\big) \big(E G E\big)\big(\frac{1}{2} D P E \big) \nonumber\\
& =  E P' G P E \nonumber\\
& =  E H E.
\mylabel{eq:6.3} 
\end{align}
That is, the equivalence class of $\Phi(g)$ does not depend on the particular choice of Convenient Shape 1. 
Therefore $\Phi$ extends to a well defined mapping from the equivalence classes of forms in $TG_{1,p}$ to forms with $\Delta = 16  p^2$ 
that are  classically integral and do not represent any numbers $n$ with $n \equiv 1,2 \pmod 4$.\\
Now, let us take the collection of all the forms with $ \Delta = 16  p^2$ that are primitive,  classically integral, and do not represent any numbers $n$ with $n \equiv 1,2 \pmod 4$ 
and call that $TG_{2,p}$. It is not difficult to show that such forms can be put into Convenient Shape 2, with Gram matrix $H$. 
Then $\frac{1}{4} D H D$ is the Gram matrix of a form in $TG_{1,p}$. 
It is not difficult to show that this "downwards" map also respects equivalence classes of forms, and the choice of Convenient Shape 2 does not matter. 
Therefore it is legitimate to name this mapping $ \Phi^{-1}$. As $\Phi$ and $\Phi^{-1}$ really are inverses, it follows that both are injective and surjective.
Suppose $g$ and $h$ are in $TG_{1,p}$  and in Convenient Shape 1, with Gram matrices $G_{1}$ and $H_{1}$, repectively. 
As they are in the same genus, there is an odd number $w$ not divisible by $p$, along with an integral matrix $R$, such that 
$$
R' G_{1} R = w^2 H_{1},
$$  
and
\beqs
\det G_{1} = \det H_{1}.
\eeqs  
This is Siegel's definition of a genus: rational equivalence ``without essential denominator''. 
Let $\Phi(g)$ have Gram matrix $G_{2}$, while  $\Phi(h)$ has Gram matrix $H_{2}$.  Then $Q =\frac{1}{2} D R E$ is integral, and we have
$$
Q' G_{2} Q = w^2 Q_{1}.
$$ 
That is, $\Phi(g)$ and $\Phi(h)$ are in the same genus, which we are calling $TG_{2,p}$.
Next, if $A$ is an automorph of $g \in TG_{1,p},$ in Convenient Shape 1, with Gram matrix $G$, this means that $A$ has determinant $\pm 1$ and 
$$
A' G A = G.
$$ 
So, repeating \eqn{6.3}, we find that
$ B = \frac{1}{2} D A E$ is an automorph of $\Phi(g)$.  At the same time, beginning with  $h \in TG_{2,p}$, in Convenient Shape 2, with Gram matrix $H$, and an automorph $B$ solving $B' H B =  H$, then
$ A = \frac{1}{2} E Q D$ is an automorph of $\Phi^{-1}(h)$. That is, the number of automorphs are the same, from which it follows that the mass of $TG_{2,p}$ is exactly the same as the mass of $TG_{1,p}$.
That is,
\beq
M(TG_{2,p})= M(TG_{1,p})=\frac{p-1}{48}.
\mylabel{eq:6.4}
\eeq
A very similar formalism shows directly that
$$
R_g(n) = R_{\Phi(g)}(4n),
$$
where these are the (finite) number of representations by the indicated form.
Now, from Theorem 29 in \cite{Wat}, we know that all forms in $TG_{1,p}$ are equivalent over the $2$-adic integers to 
$y z - x^2$, or  $\langle -1,0,0,1,0,0 \rangle$ which is integral and is in Convenient Shape 1. 
The same process that took some $g \in TG_{1,p}$ and constructed automorphs or equivalences involving $\Phi(g)$ can be readily extended to the $2$-adic integers. 
So we begin with $g \sim_2 \langle -1,0,0,1,0,0 \rangle$ which shows that 
$$
\Phi(g) \sim_2  \Phi(\langle -1,0,0,1,0,0 \rangle  ) = \langle -1,0,0,4,0,0 \rangle.
$$ 
So it follows that for any $h \in TG_{2,p}$, $h$ is equivalent over the $2$-adic integers to $4yz - x^2$, which is in Convenient Shape 2. 
So we have proved the following identities 
\begin{align}
& g\sim_{2} yz-x^2,\;\; g \in TG_{1,p},\nonumber\\
& h\sim_{2} 4yz-x^2,\;\; h \in TG_{2,p}.
\mylabel{eq:6.5}
\end{align}
Now to the $p$-adic diagonalization of these forms, which requires more terminology. 
The forms in either genus are isotropic in the $2$-adic field, as there are nontrivial integral expressions with
$yz-x^2 = 0$ or $4yz-x^2 = 0$. It follows that the forms in both genera are anisotropic (not ``zero forms'') in the $p$-adic field. 
This is from Lemma 1.1 in \cite{Cas}, page 76.
What sort of numbers are represented by these forms? 
According to Corollary 13 in \cite{Jones}, page 41, some number $n$ is represented by $g(x,y,z) \in TG_{1,p} $ in $\mathbf Q_p$ if and only if 
$$
h(x,y,z,w) = g(x,y,z) - n w^2
$$ 
is isotropic in $\mathbf Q_p$. The determinant of $h$ is $-n p^2$. This is a square in $\mathbf Q_p$ if $( -n | p) = 1$.   
We already know that $c_p(h) = c_p(g) = -1$ by Lemma 2.3(iii) in \cite{Cas}, page 58. Thus $c_p(h) =  -(-1,-1)_p$. Here $(a,b)_p$ is the Hilbert Norm Residue Symbol. 
By Lemma 2.6 on page 59, when $(-n | p )=1$ we have $h(x,y,z,w)=g(x,y,z)-n w^2$ anisotropic in ${\mathbf Q}_p$ and so $n$ is not represented.
So, for $p \equiv 1 \pmod 4$, forms in $TG_{1,p}$ and $TG_{2,p}$ represent only quadratic nonresidues modulo $p$, among the numbers not divisible by $p$. 
For $p \equiv 3 \pmod 4$, forms in $TG_{1,p}$ and $TG_{2,p}$ represent only quadratic residues modulo $p$.
Let $p$ be an odd prime and $u$ be any integer with $(-u|p)=-1$. From the fact that any binary form
$a x^2 + b x y + c y^2$ with discriminant not divisible by $p$ represents both residues and nonresidues modulo $p$, it follows that
$g \in TG_{1,p}$ diagonalizes as $u x^2 + p (v y^2 + w z^2)$. Now, given any number $V$ with
$(-V | p ) = 1$, it follows that $ x_1^2 + V x_2^2 $ is isotropic in $\mathbf Q_p$. As our forms are anisotropic there, it follows that $( -vw | p ) = -1$.  
Meanwhile, as $v,w$ are units in $\mathbf Q_p$, we know that the binary $v y^2 + w z^2$ represents both residues and nonresidues modulo $p$.
From  Lemma 3.4 in \cite{Cas}, page 115, we can insist that $v=1$ and $w = u$, so that
\beq
g\sim_p u x^2 + p ( y^2 + u z^2 ),\;\;\text{with } (-u|p)=-1,\;\; g \in TG_{1,p}.
\mylabel{eq:6.6}           
\eeq
By the usual methods, the same applies to $TG_{2,p}$. And so 
\beq
h\sim_p u x^2 + p ( y^2 + u z^2 ),\;\;\text{with } (-u|p)=-1,\;\; h \in TG_{2,p}. 
\mylabel{eq:6.7}           
\eeq

It turns out that our bijection $\Phi$ is an instance of a Watson transformation \cite{Wat2}.
As a result, the bijection generalizes to positive ternary forms with any odd discriminant.
The phenomenon of bijection of automorphs and equal mass then generalizes to any dimension and any discriminant.
We will discuss this further in Section 8.

\bigskip
\section{Proof of Theorem 1.3.} \label{sec:7}
\medskip

Here we prove our main result \eqn{1.3}. 
We recall that, thanks to Theorem 29 in \cite{Wat}, we have
\beqs
\tilde{f} \sim_{p'} x^2+y^2+z^2, 
\eeqs
for any ternary form $\tilde{f}$ with discriminant $\Delta$, provided  prime $p'\nmid 2\Delta$. \\
Hence
\beq
d_{TG_{1,p},p'}=d_{TG_{2,p},p'} =d_{x^2+y^2+z^2,p'},\;\; p'\nmid 2p.
\mylabel{eq:7.1}
\eeq
Next, we employ \eqn{2.1}, \eqn{2.6}, \eqn{3.1}, \eqn{4.3}, \eqn{6.4}, \eqn{6.5}, \eqn{6.6}, \eqn{6.7}, and \eqn{7.1} to rewrite the expression on the right of \eqn{1.3} as
\begin{align*}
&\mbox{RHS(1.3)}=\\
& \frac{p-1}{48} 96\pi \sqrt{\frac{n}{p^2}} (2 d_{TG_{1,p},2}(n) - d_{TG_{2,p},2}(n)) \frac{p}{p-1}\Gamma_{p}(n)\prod_{\gcd(p',2p)=1}d_{x^2+y^2+z^2,p'}(n) \\
& = 2\pi \sqrt{n} (2 d_{TG_{1,p},2}(n) - d_{TG_{2,p},2}(n)) \Gamma_{p}(n)\prod_{\gcd(p',2p)=1}d_{x^2+y^2+z^2,p'}(n)\\
& = 2\pi \sqrt{n} \psi(n) \Gamma_{p}(n)\prod_{\gcd(p',2p)=1}d_{x^2+y^2+z^2,p'}(n) =\mbox{LHS(1.3)}.
\end{align*}
This completes our proof of the Theorem 1.3.\\
We conlude with the following example.
Genus $TG_{1,73}$ consists of four classes
\beqs
TG_{1,73} = \{\mbox{Cl}(h_{1}), \mbox{Cl}(h_{2}), \mbox{Cl}(h_{3}),\mbox{Cl}(h_{4})\},
\eeqs
where 
\beqs
\begin{array}{ll}
h_{1}(x,y,z) = 31x^2+5y^2+11z^2 + yz    -14zx + 6 xy, |\mbox{Aut}(h_{1})|=2,\\
h_{2}(x,y,z) = 15x^2+14y^2+10z^2 + 7yz  +4zx + 16xy, |\mbox{Aut}(h_{2})|=2,\\
h_{3}(x,y,z) = 11x^2+7y^2+20z^2 + 7yz   +2zx  +  4xy, |\mbox{Aut}(h_{3})|=4,\\
h_{4}(x,y,z) = 7x^2+11y^2+21z^2 + 11yz   +2zx  +  4xy, |\mbox{Aut}(h_{4})|=4.\\
\end{array}
\eeqs
Note that all four forms above are in Convenient Shape 1.
And so, we can immediately construct the second genus
\beqs
TG_{2,73} = \{\mbox{Cl}(g_{1}), \mbox{Cl}(g_{2}), \mbox{Cl}(g_{3}),\mbox{Cl}(g_{4})\},
\eeqs
where 
\beqs
\begin{array}{ll}
g_{1}(x,y,z) = 31x^2+20y^2+44z^2 + 4yz    -28zx + 12 xy, |\mbox{Aut}(g_{1})|=2,\\
g_{2}(x,y,z) = 15x^2+56y^2+40z^2 + 28yz  +8zx + 32xy, |\mbox{Aut}(q_{2})|=2,\\
g_{3}(x,y,z) = 11x^2+28y^2+80z^2 + 28yz   +4zx  +  8xy, |\mbox{Aut}(g_{3})|=4,\\
g_{4}(x,y,z) = 7x^2+44y^2+84z^2 + 44yz   +4zx  +  8xy, |\mbox{Aut}(h_{4})|=4.\\
\end{array}
\eeqs
From \eqn{1.3} we get 
\begin{align}
s(73^2 n)-73s(n)& = 24(31,5,11,1,-14,6)(n)+24(15,14,10,7,4,16)(n) \nonumber\\ 
& + 12(11,7,20,7,2,4)(n) + 12(7,11,21,11,2,4)(n) \nonumber\\
& - 48(31,20,44,4,-28,12)(n) - 48(15,56,40,28,8,32)(n) \\
& - 24(11,28,80,28,4,8)(n) - 24(7,44,84,44,4,8)(n).\nonumber
\mylabel{eq:7.2}
\end{align}

\bigskip
\section{Watson's Transformations}\label{sec:8}
\medskip

 We begin with a brief summary of the Watson ``$m$-mapping''.

Let $f$ be an $n$-ary quadratic form. In \cite{Wat2} equation (2.4), Watson defines a certain lattice $\Lambda_m(f)$ by $ \vec{x} \in\Lambda_m(f)$ if and only if
$ f(\vec{x} + \vec{z}) \equiv  f(\vec{z}) \pmod m, \; \; \forall \vec{z} .$ This is a lattice, because if $\vec{x} \in \Lambda_m(f),$ then $-\vec{x} \in \Lambda_m(f),$ and if $\vec{x_1},\vec{x_2} \in \Lambda_m(f),$ then $\vec{x_1}+ \vec{x_2} \in \Lambda_m(f)$.
Let $F$ be the Gram (or second partials) matrix for the $n$-ary form $f$. This means, with column vector $\vec{x}$ and its transposed row vector $\vec{x}'$, that 
$ f(\vec{x}) = \frac{1}{2} \; \vec{x}'  F \vec{x}$.
Watson also supplies his (2.2), which is simply
$ f(\vec{x} + \vec{z}) = f(\vec{x}) +  \vec{z}'  F \vec{x}  + f(\vec{z})$.
Now, we can vary $\vec{z}$ arbitrarily in $ f(\vec{x}) +  \vec{z}'  F \vec{x}  + f(\vec{z})\equiv f(\vec{z}) \pmod m$. 
So $ \vec{x} \in\Lambda_m(f)$
if and only if $f(\vec{x}) \equiv 0 \pmod m$ and $ F \vec{x} \equiv \vec{0} \pmod m,$ thus giving $n+1$ equations $\pmod m$, the first one tending to be redundant 
(depending on powers of $2$ in $m$ and various coefficients).
Then Watson says to take a square integral matrix $M$ whose columns serve as an integral basis for
$\Lambda_m(f),$
and in (2.5) defines the result of the ``$m$-mapping'', that we denote  $g = \lambda_m(f)$, by
$ g(\vec{y}) = \frac{1}{m} \; f( M \vec{y})$.
So, if we take $G$ as the Gram matrix for $g$, we have
$G = \frac{1}{m} M' F M$,
where $M'$ refers to the transpose of $M$. By construction, $G$ is an integral matrix, as all entries of $FM$ are divisible by $m$,
and $ g(\vec{y}) = \frac{1}{2} \; \vec{y}'  G \vec{y}$. Watson shows that a different choice of basis matrix $M$ simply gives a form equivalent to $g$. So this is Watson's ``$m$-mapping''.

We are ready to show that  our bijection $\Phi$ is Watson's $\lambda_4$.
Let our discriminant $\Delta = \delta$ be odd, and let our form $\tilde{f}$ be in the analogue of convenient shape $1$. 
That is,  $a,d$ odd and $e,f$ even in $\langle a,b,c,d,e,f \rangle$. Notice that 
$$ 
\delta \equiv - a d^2 \equiv -a \pmod 4.
$$ 
For the case of  $TG_{1,p}$ we had $  \delta \equiv 1 \pmod 4$ so it followed that, for  $TG_{1,p}$, in convenient shape $1$ we had  $ a \equiv 3 \pmod 4$. 
But we will need only odd $a$. We just use $a\equiv -\delta  \pmod 4$.
We find that $ \vec{x} \in \Lambda_4(\tilde{f})$ when $ \tilde{f}(\vec{x}) \equiv  0 \pmod 4$  and
$$  
F \vec{x} \equiv  \vec{0} \pmod 4.
$$ 
Here  $\vec{x}$ is the column vector
\[ \vec{x} \; = \;
\left(  \begin{array}{c}
  x \\
  y\\
  z  
\end{array} 
  \right) .
  \]
So we have three linear equations $\pmod 4$. From $a,d$ odd, $e,f$ even, and $ f x + 2 b y + d z \equiv 0 \pmod 4 $ we find that $z$ is even.
From $ e x + d y + 2 c z \equiv 0 \pmod 4$ we find that $y$ is also even.
With $y,z$ even,  $a,d$ odd, $e,f$ even, and $2 a x + f y + e z \equiv 0 \pmod 4$ we find that $x$ is even.
With $x,y,z$ even,  $a,d$ odd, $e,f$ even, back to $f x + 2 b y + d z \equiv 0 \pmod 4$, we find that $z \equiv 0 \pmod 4$.
With $x,y,z$ even,  $a,d$ odd, $e,f$ even, back to $ e x + d y + 2 c z \equiv 0 \pmod 4$, we find that $y \equiv 0 \pmod 4$. 
Finally, with $x,y,z$ all even, we do indeed have $ \tilde{f}(\vec{x}) \equiv  0 \pmod 4$.

Now, look at Watson's (2.5), his $M$ must have columns that form an integral basis of $\Lambda_4(\tilde{f}),$ and this particular time it makes sense to  choose a matrix already in Smith Normal Form,
\begin{equation}
M \, = \,
\left(  \begin{array}{ccc}
  2 & 0 & 0\\
  0 & 4 & 0\\
  0 & 0 & 4  
\end{array} 
  \right). 
\end{equation}
Back to  Watson's (2.5), we get the result of the ``$m$-mapping'', that we write as $g = \lambda_4(\tilde{f})$, 
$ g(\vec{y}) = \frac{1}{4} \; \tilde{f}( M \vec{y})$, or, reverting to vector $(x,y,z)$
\begin{equation}
\label{x2y2z}
g(x,y,z) = \frac{1}{4} \; \tilde{f}(2x, 4y,4z) = \tilde{f}(x,2y,2z).
\end{equation}
The coefficients for $g$ are $ \langle a,4b,4c,4d,2e,2f \rangle$ so that all but the first are divisible by $4$, and the form $g$ cannot take any value $2,\delta \pmod 4,$ as  $a \equiv -\delta  \pmod 4$. 
Meanwhile (\ref{x2y2z}) shows that $g$ represents a subset of the numbers represented by $\tilde{f}$.
Also (\ref{x2y2z}) shows that, in the case of $TG_{1,p},$ $\Phi = \lambda_4$.

Our new discriminant is $\Delta = 16 \delta$, with any form $g$ in convenient shape $2$, that is
$\langle a,b,c,d,e,f \rangle  $ with some odd  $ a  \equiv -\delta  \pmod 4$ and $b,c,d,e,f$ divisible by $4$, so that $g \neq 2,\delta \pmod 4$. 
Note that, with odd $\delta$, and $\Delta \equiv \pm 16 \pmod{64}$, we know that this $d \equiv  4 \pmod 8$.
So, if we call this form $g(x,y,z)$, we have $g(x,y,z) \equiv a x^2 \pmod 4$. 
If we let $\vec{z}$ be the column vector with entries $u,v,w$, Watson's (2.4) for  $\Lambda_4(g)$ becomes $a (x+u)^2 \equiv a u^2 \pmod 4, \; \forall u$. 
So, $ a x^2 + 2 a u x \equiv 0 \pmod 4, \; \forall u$, and it is necessary and sufficient to have the first coordinate $x$ even. So this time, using the matrix name $N$, we get
\begin{equation}
N \, = \,
\left(  \begin{array}{ccc}
  2 & 0 & 0\\
  0 & 1 & 0\\
  0 & 0 & 1  
\end{array} 
  \right), 
\end{equation}
and $MN=NM=4I.$

Let us take, say, $h = \lambda_4(g)$.
Back to  Watson's (2.5), we get $ h(\vec{y}) = \frac{1}{4} \; g( N \vec{y})$, or, reverting to vector $(x,y,z)$
$$ 
h(x,y,z) = \frac{1}{4} \; g(2x, y,z) = g(x,\frac{y}{2},\frac{z}{2}).
$$
So, once again, in the case of $TG_{1,p},$ $\Phi^{-1} = \lambda_4.$ Furthermore, beginning with a genus of any odd discriminant $\delta,$ $\lambda_4^2(\tilde{f}) = \tilde{f}$.

  We now proceed to describe a very general result that could easily have been present in \cite{Wat2} but is not there in any explicit manner, and indeed for which we have no reference. So this result about automorphs and mass may be new to some extent. Note that Watson uses $(f)$ to refer to the integral equivalence class of the form $f$.
\noindent
\begin{theorem}\label{nt8}
In any dimension $n,$ with any $m$ and any discriminant, if $\lambda_m^2$ is the identity on a class $(f),$ then $\lambda_m$ is a bijection between the sets of equivalence classes in the genus of $f$ and the genus of $\lambda_m(f)$. Furthermore, $\lambda_m$ induces a bijection of integral automorphs between any form $f_1$ in the genus of $f$ and its image $\lambda_m(f_1)$. 
As a result, taking positive forms, the two genera have the same mass.
\end{theorem}
\noindent
Our proof is long but elementary, and can be written entirely using Watson's terminology and concepts. We content ourselves with a sketch of the proof.
Let the columns of the matrix $M$ be an integral basis for $\Lambda_m(f)$, and let $g = \lambda_m(f)$, so that if $G$ as the Gram matrix for $g$, we have $G = \frac{1}{m} M' F M$.  
Watson  shows that equivalent forms are mapped to equivalent forms, forms in the same genus are mapped to forms in the same genus, and that the map from genus to genus is surjective. It is not usually injective, and usually $\lambda_m^2\left( (f) \right)$ has a different discriminant from that of $f$, therefore usually can not even be in the same genus as $f$.
If it should happen that $\lambda_m^2\left( (f) \right) = (f)$, there is, in fact, a matrix $N$, where $N$ turns out to be integral as well (Watson's Lemma 2(ii), page 581), such that
\begin{eqnarray}
 MN=NM=mI, \\
 G = \frac{1}{m} M' F M, \\
 F = \frac{1}{m} N' G N.
\end{eqnarray}
These formulae extend to a bijection between genera, where the actual matrix $M$ depends on $f$, and of course $N$ depends on $M$.
Finally, if $\lambda_m(f) = g$, and if $R$ is an integral automorph of $f$,  that is $R' F R = F$, then
\begin{equation}
S = \frac{1}{m} N R M 
\end{equation} 
is an integral automorph of $g$, that is $S' G S = G$. And if $S$ is an integral automorph of $g$, then
\begin{equation}
R = \frac{1}{m} M S N 
\end{equation} 
is an integral automorph of $f$.
The point that required extra care in the proof was the fact that, in our situation, the columns of the matrix $N$ really do give an integral basis for $\Lambda_m(g)$.
The key to this is  Watson's sentence on page 585, 
``Now the foregoing argument shows that we have $|H| = \pm 1$ and $f\sim \phi$ if and only if the lattice $ \mu_m(f)$ is $m\Lambda_n$''.This refers to his Theorem 2 on page 580.

\bigskip
\noindent
\textbf{Acknowledgements}
\medskip

\noindent
We are grateful to Billy Chan, Larry Gerstein, and Ken Ono for their interest and comments.
We would like to thank Rainer Schulze-Pillot for bringing \cite{Leh} to our attention and Keith Grizzell for careful reading of the manuscript.

\bigskip

\bibliographystyle{amsplain}

\begin{thebibliography}{16}
%
\bibitem{Bach}
P. Bachmann, \emph{Die Arithmetik von Quadratischen Formen}, Leipzig, 1898.
%
\bibitem{PB}
P.T. Bateman, \emph{On the representations of a number as the sum of three squares}, Trans. Amer. Math. Soc. \textbf{71} (1951), 70--101.
%
\bibitem{Brk}
A. Berkovich, \emph{On representation of an integer by $X^2+Y^2+Z^2$ and the modular equations of degree $3$ and $5$}, arXiv:0907.1725v3,
to appear in the volume "Quadratic and Higher Degree Forms", in Developments in Math., Springer 2011.
%
\bibitem{BeJ}
A. Berkovich, W.C. Jagy,
\emph{Ternary quadratic forms, modular equations and certain positivity conjectures}, in K. Alladi, J. Klauder, and C.R. Rao, eds, 
The Legacy of Alladi Ramanakrishnan in the mathematical sciences, 211--241, Springer, NY, 2010.
%
\bibitem{BHJ}
A. Berkovich, J. Hanke, W.C. Jagy,
\emph{A proof of the $S$-genus identities for ternary quadratic forms}, arXiv:1010.1926.
%
\bibitem{Cas}
J.W.S. Cassels, \emph{Rational Quadratic Forms}, Dover, 2008.
%
\bibitem{CS}
J.H. Conway, N.J.A. Sloane,
\emph{Low-dimensional lattices.IV. The mass formula}, Proc. Roy. Soc. London Ser. A \textbf{419} (1988), 259--286.
%
\bibitem{Dick}
L.E. Dickson, \emph{Modern Elementary Theory of Numbers}, The University of Chicago Press, 1939.
%
\bibitem{Jones}
B.W. Jones, \emph{The Arithmetic Theory of Quadratic Forms}, Mathematical Association of America, 1950.
%
\bibitem{Landau}
E.Landau, \emph{Vorlesungen \"uber Zahlentheorie}, Leipzig: S. Hirzel, 1927.
%
\bibitem{Leh}
J.L. Lehman, \emph{Levels of positive definite ternary quadratic forms}, Math. of Comput. \textbf{58} (1992), 399--417.
%
\bibitem{Siegel} 
C.L. Siegel, \emph{Lectures on the Analytical Theory of Quadratic Forms}, Notes by Morgan Ward.
Third revised edition. Buchhandlung Robert Peppm\"uller, G\"ottingen, 1963.
%
\bibitem{Wat} 
G.L. Watson, \emph{Integral Quadratic Forms}, Cambridge University Press, 1970.
%
\bibitem{Wat2} 
G.L. Watson, \emph{Transformations of a quadratic form which do not increase the class-number}, Proc. London Math. Soc. \textbf{12} (1962), 577--587.

http://plms.oxfordjournals.org/content/s3-12/1/577.full.pdf

\end{thebibliography}

\end{document}